\documentclass[final,onefignum,onetabnum]{siamart171218}

\usepackage{lipsum}
\usepackage{amsfonts}
\usepackage{graphicx}
\usepackage{epstopdf}
\usepackage{algorithmic}
\usepackage{booktabs}
\usepackage{pgfplots,pgfplotstable}
\pgfplotsset{compat=1.18}
\usepackage{listings}
\usepackage{multirow}
\usepackage[nolist,nohyperlinks]{acronym}

\newcommand{\revNB}[1]{\textcolor{red}{#1}}

\acrodef{dof}[DoF]{Degree of Freedom}
\acrodefplural{dof}[DoFs]{Degrees of Freedom}
\acrodef{aap}[AAP]{Alternating Anderson-Picard}
\acrodef{aar}[AAR]{Alternating Anderson-Richardson}
\acrodef{aa}[AA]{Anderson acceleration}
\acrodef{pde}[PDE]{Partial Differential Equation}
\acrodefplural{pde}[PDEs]{Partial Differential Equations}

\lstdefinelanguage{Julia}%
  {morekeywords={abstract,break,case,catch,const,continue,do,else,elseif,%
      end,export,false,for,function,immutable,import,importall,if,in,%
      macro,module,otherwise,quote,return,switch,true,try,type,typealias,%
      using,while},%
   sensitive=true,%
   alsoother={$},%
   morecomment=[l]\#,%
   morecomment=[n]{\#=}{=\#},%
   morestring=[s]{"}{"},%
   morestring=[m]{'}{'},%
   backgroundcolor = \color{blue!4!white},
}[keywords,comments,strings]%

\ifpdf
  \DeclareGraphicsExtensions{.eps,.pdf,.png,.jpg}
\else
  \DeclareGraphicsExtensions{.eps}
\fi

\newsiamremark{remark}{Remark}
\newsiamremark{hypothesis}{Hypothesis}
\crefname{hypothesis}{Hypothesis}{Hypotheses}
\newsiamthm{claim}{Claim}

\headers{Anderson acceleration for coupled physics applications}{Nicolas A. Barnafi, Massimiliano Lupo Pasini}

\title{Convergence analysis and considerations on computational efficiency of physics-aware Alternating Anderson acceleration for coupled physics applications\thanks{Submitted to the editors DATE.
}}

\author{Nicolas A. Barnafi\thanks{Imagination Corp., Chicago, IL 
  (\email{nicolas.barnafi@uc.cl}, \url{http://www.imag.com/\string~ddoe/}).}
\and Massimiliano Lupo Pasini\thanks{Department of Applied Mathematics, Fictional University, Boise, ID 
  (\email{lupopasinim@ornl.gov}.}
}

\usepackage{amsopn}

\DeclareMathOperator{\argmin}{argmin}

\makeatletter
\newcommand*{\addFileDependency}[1]{
  \typeout{(#1)}
  \@addtofilelist{#1}
  \IfFileExists{#1}{}{\typeout{No file #1.}}
}
\makeatother

\newcommand*{\myexternaldocument}[1]{%
    \externaldocument{#1}%
    \addFileDependency{#1.tex}%
    \addFileDependency{#1.aux}%
}

\renewcommand{\vec}{\vectorsym}
\newcommand{\mat}{\matrixsym}
\newcommand{\ten}{\tensorsym}
\newcommand{\parder}[2]{\frac{\partial #1}{\partial {#2}}}

\DeclareMathOperator{\grad}{\nabla}
\DeclareMathOperator{\dive}{\text{div}}

\newcommand{\R}{\mathbb{R}}

\renewcommand{\vec}{\vectorsym}
\newcommand{\pressureS}{p_{S}}
\newcommand{\pressureNS}{p_{NS}}
\newcommand{\pL}{p_{L}}

\ifpdf
\hypersetup{
  pdftitle={Anderson acceleration for coupled physics applications},
  pdfauthor={Nicolás A. Barnafi, Massimiliano Lupo Pasini}
}
\fi

\myexternaldocument{ex_supplement}

\usepackage{graphicx} 
\usepackage{amsmath, amssymb, mathpazo, isomath, mathtools}
\usepackage{subcaption,graphicx,pgfplots}
\usepackage{booktabs}
\usepackage{hyperref}

\title{Two-Level Sketching Alternating Anderson acceleration for Complex Physics Applications}
\author{Nicolás A. Barnafi\thanks{Instituto de Ingeniería Matemática y Computacional \& Facultad de Ciencias Biológicas, Pontificia Universidad Católica de Chile, Santiago, RM Chile, and Centro de Modelamiento Matemático, Santiago, RM Chile.} \and Massimiliano Lupo Pasini\thanks{Computational Sciences and Engineering Division, Oak Ridge National Laboratory, Oak Ridge, TN 37831  USA.}}
\date{}

\begin{document}

\maketitle

\begin{abstract}
We present a novel two‐level sketching extension of the Alternating Anderson–Picard (AAP) method for accelerating fixed‐point iterations in challenging single‐ and multi‐physics simulations governed by discretized partial differential equations. Our approach combines a static, physics‐based projection that reduces the least‐squares problem to the most informative field (e.g., via Schur–complement insight) with a dynamic, algebraic sketching stage driven by a backward stability analysis under Lipschitz continuity. We introduce inexpensive estimators for stability thresholds and cache‐aware randomized selection strategies to balance computational cost against memory‐access overhead. The resulting algorithm solves reduced least‐squares systems in place, minimizes memory footprints, and seamlessly alternates between low‐cost Picard updates and Anderson mixing. Implemented in Julia, our two‐level sketching AAP achieves up to 50\% time‐to‐solution reductions compared to standard Anderson acceleration—without degrading convergence rates—on benchmark problems including Stokes, $p$‐Laplacian, Bidomain, and Navier–Stokes formulations at varying problem sizes. These results demonstrate the method’s robustness, scalability, and potential for integration into high‐performance scientific computing frameworks.
\end{abstract}

{\footnotesize \noindent This manuscript has been authored in part by UT-Battelle, LLC, under contract DE-AC05-00OR22725 with the US Department of Energy (DOE). The US government retains and the publisher, by accepting the article for publication, acknowledges that the US government retains a nonexclusive, paid-up, irrevocable, worldwide license to publish or reproduce the published form of this manuscript, or allow others to do so, for US government purposes. DOE will provide public access to these results of federally sponsored research in accordance with the DOE Public Access Plan (\url{http://energy.gov/downloads/doe-public-access-plan}).}

\section{Introduction}
The \ac{aa} method \cite{Anderson} is a multi-secant method \cite{brezinski_shanks_2018, saad2021} that has been widely used either to improve the convergence rate of convergent fixed-point schemes or to restore convergence when the original fixed-point scheme is not convergent \cite{Fang, walker_anderson_2011, potra_characterization_2013, zhang2020, saad2021}. In particular, the convergence of \ac{aa} has been studied with respect to specific properties of different scientific applications \cite{toth_convergence_2015, toth_local_2017, Mai2019, doi:10.1137/20M132938X}.

\ac{aa} requires solving a least-squares (LS) {problem} at each fixed-point iteration where up to all previous iterations are used. This can be computationally expensive for scientific applications that involve large-scale calculations, especially when such calculations are distributed on high-performance computing (HPC) platforms.
To reduce the computational cost of \ac{aa} for the solution of linear systems, a recently proposed variant of \ac{aa} called \ac{aar} \cite{banerjee_periodic_2016, AAJ, AAR} performs multiple (inexpensive) Richardson iterations between two consecutive Anderson mixing corrections. This can be naturally extended to the non-linear setting by using Picard iterations, and is known as the \ac{aap} method \cite{laiu2024aap}. The number of fixed-point iterations between two consecutive Anderson mixing corrections can be arbitrarily set by the user. When no such iterations are performed, \ac{aap} becomes \ac{aa}. Recent theoretical results establish that for linear systems \ac{aar} is more robust than the standard \ac{aa} against stagnations \cite{lupo_pasini_convergence_2019}.  However, even if performed intermittently, computing the Anderson mixing via least-squares can still be computationally expensive when the number of rows in the tall and skinny matrix that defines the left-hand side (LHS) of the LS problem is large. 

Various literature contributions assessed the computational advantages of replacing accurate and expensive computations with inaccurate and inexpensive ones in \ac{aap}, either by approximating the evaluation of the fixed-point operator or by reducing the accuracy to solve the LS problem. 
For specific problems, the accuracy and computational costs are controlled using physics information, which results in computational savings without affecting the final convergence of the fixed-point scheme \cite{doi:10.1137/19M1290097, doi.org/10.1111/cgf.14081}. 
For situations where physics-driven guidelines are difficult (and often impossible) to determine, \cite{lupo2022anderson} proposed more general guidelines motivated by forward and backward error analysis, which propose to  judiciously subselect the rows of the tall and skinny matrix associated with the LS problem via adaptive sketching \cite{woodruff2014sketching, mahoney2011randomized, murray2023randomized}. 

In this work, we extend the backward stability analysis performed in \cite{lupo2022anderson} to non-linear fixed-point iterations under the assumption of Lipschitz continuity of the fixed-point operator, and we propose a two-level sketching technique to project the LS problem in a subspace and thereby reduce the computational cost of solving it while still ensuring convergence of \ac{aap} to the desired accuracy. The first sketching is motivated by the physics, and consists in projecting the residual onto a subset of the fields that characterize the physics of the phenomenon of interest. The use of the first sketching level is motivated by the fact that (i) multi-physics problems can be recast as a system of only one of their variables without losing information (such as the Schur complement in linear systems), and (ii) different fields can attain values with very different scales, so that combining them can hinder the overall conditioning of the LS problems solved to compute the Anderson mixing, which in turn can affect the convergence of \ac{aap}. The first level physics-based sketching operator is static and determines a fixed projection subspace throughout the solution process, which allows to conveniently pre-allocate the memory needed once at the beginning of the execution of the algorithm. The second sketching level of our approach coincides with the adaptive sketching already used in \cite{lupo2022anderson}, for which we propose sharper indicators to activate an adaptive selection of the rows of the tall an\revNB{d} skinny matrix according to our backward stability analysis.

While each of the two components—projection of Anderson acceleration onto a reduced or physically meaningful subspace and approximate or sketched solution of the associated least-squares problem—has appeared separately in prior work (e.g., projection and variable-aware formulations of Anderson acceleration \cite{walker_anderson_2011, 10.1093/imanum/draa095} and backward-stable inexact or perturbed Anderson least-squares solves \cite{lupo2022anderson}), the contribution of this paper is not merely their juxtaposition. Rather, the manuscript introduces a structured two-level decomposition of approximation error in Anderson acceleration, in which the static physics-based restriction and the dynamic algebraic sketching address distinct and complementary stability mechanisms. The backward stability analysis clarifies that physics-based projections act on the conditioning and scale of the least-squares problem, while algebraic sketching introduces controlled perturbations whose admissibility depends on residual-dependent bounds. This separation leads to an algorithmic structure in which the two sketching levels play fundamentally different roles and cannot be interchanged or merged without loss of robustness. We highlight that the resulting adaptive scheme can still suffer from instabilities that are inherent to the method, as the backward stability analysis only ensures that the sketching does not result in \emph{additional} instability. Other improvements have been considered in the literature, such as better memory management in the QR factorization through a Truncated Gram-Schmidt approach \cite{tang2024anderson}, as well as enriching the method with safeguarding mechanisms that provide a strict improvement on convergence \cite{garstka2022safeguarded}.

When the sketching at the second level is performed by randomly selecting the entries to retain, there is an increased risk that the entries retained are physically stored far from each other in the computer memory, which in turn increases the chances of cache misses and slows down the execution of the code. To address this issue, we propose simple but effective guidelines for the randomized sketching that reduce the overhead of cache-misses in two specific operations: matrix-vector products and the QR factorization. To define a sketching range, we compare the time elapsed to perform each operation using sequential memory access against a random mask, and from there obtain the maximum sketching allowed to have faster operations for various problem sizes. Based on benchmark tests that we performed in this work, we concluded that a reasonable percentage of entries to retain in the second (adaptive) level sketching is 10\% for matrix-vector products and 80\% for QR factorizations. While the arithmetic complexity of the least-squares solve scales linearly with problem size, the practical performance is dominated by memory traffic and cache coherency due to the non-contiguous memory access patterns of the masks. Consequently, in this work we rely on empirical benchmarks rather than asymptotic analysis to characterize the method's efficiency.

We implemented our approach in the coding language \texttt{Julia} \cite{Bezanson_Julia_A_fresh_2017}, and illustrated the performance of our numerical scheme on different single- and multi-physics problems: Stokes, $p$-Laplacian, Bidomain, and stationary Navier-Stokes. We additionally considered an increasing number of \acp{dof} in all test, in order to ensure that the proposed methods do not present deterioration for increasing problem sizes. Our numerical results confirmed the efficacy of our sketching approach in both reducing memory usage and computational time in \ac{aap} while still maintaining its convergence to the desired accuracy. {Our implementation is open-source and available at the following GitHub repository \url{https://github.com/ORNL/Restricted-Alternating-Anderson-Picard}}.

{From a methodological perspective, this work contributes a stability-driven framework for sketching Anderson acceleration, in which backward error bounds for the Anderson least-squares problem are used to derive admissible and adaptive row-selection strategies for nonlinear fixed-point iterations.}

The remainder of the paper is structured as follows. In Section~\ref{section:aar} we review \ac{aap}. In Section~\ref{section:aar sketching}, we motivate and propose our two-level sketching algorithm for \ac{aap}, together with the relevant notation and memory requirements. In Section~\ref{section:sketching analysis} we perform the backward stability analysis of the LS problem solved to computed the Anderson mixing vector in \ac{aap}. Based on the conclusions of the backward stability analysis we propose two sketching strategies, subselection and randomized, which extend the ones proposed in \cite{lupo2022anderson} to non-linear fixed point problems under the assumption that the fixed-point operator is Lipschitz. Finally, we provide several numerical tests in Section~\ref{section:tests} to validate our proposed methodology and briefly conclude our work with possible future research directions in Section~\ref{section:conclusions}.

\section{Alternating Anderson-Picard method}\label{section:aar}
Consider a non-linear function $\vec{T}:\R^n \to \R^n$ that define the  non-linear system
\begin{equation}
    \vec T(x) = \vec 0.
    \label{sec2:fixedpoint}
\end{equation}

We consider the Picard iteration associated with \eqref{sec2:fixedpoint}
    \begin{equation}
    \vec x^{k+1} = \vec x^k - \omega \vec T(\vec x^k) \coloneqq \vec g(\vec x^k),
    \label{sec2:picard}
    \end{equation}
    where $\omega>0$ is a relaxation parameter, so that solving the non-linear system is equivalent to finding a fixed-point of $\vec g$. 
    Denote the current fixed point residual as
        \begin{equation}
        \vec f^k \coloneqq T(\vec x^k)  
    \end{equation}
    and the current fixed point iterate as
    \begin{equation}
    \vec g^k \coloneqq \vec x^k - \omega \vec f^k.
    \end{equation}

    We then define the increment matrices

    \begin{align}
        \mat X^k&\coloneqq[\Delta \vec x^{k-m_k} \hdots \Delta \vec x^{k-1}] \in \mathbb{R}^{n\times m_k}, \\
        \mat F^k &\coloneqq [\Delta \vec f^{k-m_k} \hdots \Delta \vec f^{k-1}] \in \R^{n\times m_k},  \\
        \mat G^k &\coloneqq [\Delta \vec g^{k-m_k} \hdots \Delta \vec g^{k-1}] \in \R^{n\times m_k}, 
    \end{align}
    where  $\Delta (\cdot)^k \coloneqq (\cdot)^k - (\cdot)^{k-1}$,  $m_k \coloneqq \min\{m,k\}$ and $m$ is an integer that describes the maximum number of previous terms of the sequence. Consider the LS problem 
    \begin{equation}\label{eq:ls-aar}
        \min_{\vec\alpha}|\mat F^k\vec \alpha - \vec f^k|,
    \end{equation}
    whose solution we will denote by 
    \begin{equation} \vec \alpha^k \coloneqq \left[\mat F^k\right]^\dag \vec f^k,\end{equation}
    where $[\mat F^k]^\dagger\coloneqq ([\mat F^k]^T[\mat F^k])^{-1}[\mat F^k]^T$ is the Moore-Penrose pseudo-inverse \cite{golub2013matrix}. The solution to the LS problem $\vec\alpha$ contains the weights for the Anderson mixing, which is performed as follows:
    \begin{equation}
    \vec x^{k+1} = \vec x^k - \omega\vec f^k + \mat G^k\vec \alpha^k.
    \end{equation}
    Algorithm~\ref{alg:aar} describes all the steps performed by the \ac{aap} method, where the parameter $p$ represents the number of Picard iterations separating two consecutive \ac{aa} steps.

    \begin{algorithm}
        \caption{Alternating Anderson-Picard method (AAP)}\label{alg:aar}
        \begin{algorithmic}[1]
            \STATE{Given $\vec x^0$, compute $\vec f^0=\vec T(\vec x^0)$ and $\vec x^1=\vec x^0-\vec f^0$}
            \STATE{Set $\texttt{error}=1$, $\texttt{it}=0$}
            \WHILE{$\texttt{error} \geq \texttt{tol}$ and $\texttt{it} \leq \texttt{maxit}$}
                \IF{$k \neq 0 \mod p$}
                    \STATE{Evaluate residual $\vec f^k = \vec T(\vec x^k)$} 
                    \STATE{Update $\vec x^{k+1} = \vec x^k - \omega \vec f^k$}
                \ELSE
                    \STATE{Evaluate residual $\vec f^k = \vec T(\vec x^k)$} 
                    \STATE{Compute least squares solution $\vec \alpha^k = \left[\mat F^k\right]^\dag \vec f^k$}
                    \STATE{Update $\vec x^{k+1} \gets \vec x^k - \omega\vec f^k + \mat G^k\vec \alpha^k$}
                \ENDIF
                \STATE{Update \texttt{error} and $\texttt{it}$}
            \ENDWHILE
            \RETURN $\vec x^*$
        \end{algorithmic} 
    \end{algorithm}
    If the Anderson mixing is performed at each iteration (i.e., $p=1$), then \ac{aap} becomes the standard \ac{aa}. If $g$ is linear and $p=1$, then the fixed-point scheme becomes the Richardson's method to solve linear systems, and \ac{aa} is equivalent to GMRES in exact arithmetic and without restart, as long as there are no stagnations in the residual across two consecutive iterations of the fixed-point scheme \cite{walker_anderson_2011}. If $g$ is linear and $p>1$, then \ac{aap} becomes \ac{aar}, which is more robust against stagnations in the residual than the traditional \ac{aa} \cite{lupo_pasini_convergence_2019}, and is still equivalent with GMRES at periodic intervals, that is every time the Anderson mixing is performed after $p$ Picard steps. 
    If $\vec g$ is non-linear, then \ac{aap} with $p=1$ is equivalent to a multisecant-GMRES method that uses GMRES to
    solve a multisecant linear system at each iteration \cite{walker_anderson_2011}. If $p>1$, then the equivalence is restored every time \ac{aap} performs and Anderson update \cite{laiu2024aap}. Previous work has shown that the convergence rate of \ac{aap} with $p=1$ is at most quadratic \cite{evans2020proof}, which motivates follow-up works to use \ac{aa} only to accelerate Picard iterations with linear convergence rate \cite{rebholz2021iterated, pollock2021natural}. Moreover, \ac{aap} can increase the radius of convergence of the fixed-point iteration, and restore convergence when the underlying fixed point scheme is not convergent due to the fixed-point operator not being contractive \cite{lupo_pasini_convergence_2019, lupo2022anderson}. One commonly acknowledged feature of \ac{aap} is that alternation not only reduces the burden of solving least-squares problems, but additionally provides in many cases a \emph{reduction} of the total number of iterations \cite{banerjee_periodic_2016,AAR, lupo_pasini_convergence_2019}, which we also show in our numerical tests. A thorough understanding of what causes this improved convergence behavior for $p>1$ is still object of active research.

    Although the algorithm is relatively simple, its implementation can still significantly benefit from several hardware-aware optimizations to improve its execution. In this work, we propose numerical techniques that aim at simultaneously maximizing memory reuse and minimizing the computational cost of each \ac{aa} update. To this aim, we note that Algorithm~\ref{alg:aar} uses dense matrices, namely $\mat F^k$, $\mat X^k$, and $\mat G^k$, which are updated at each iteration, but only $\mat F^k$ and $\mat G^k$ are used for computing the Anderson mixing vector $\vec \alpha^k$. The use of sketching operators for the LS problem offers the possibility to solve the LS problem on a projected subspace by retaining only a subset of the rows, which allows to allocate only a part of the residuals matrix $\mat F^k$, and separately store the matrix $\mat G^k \coloneqq \mat X^k - \omega\mat F^k$. This opportunity of decoupling the memory storage of $\mat F^k$ and $\mat G^k$ motivates our two-level sketching approach, which allows us to allocate a smaller matrix $\mat F_\Pi$ for the LS problem, while keeping the full $\mat G^k$ matrix for the update steps. Our implementation will also perform in-place operations to further reduce the memory usage for all matrix-vector and vector-vector operations. 

\begin{remark}
    The mathematical formulation of \ac{aap} that we adopted so far follows the residual formulation described in \eqref{sec2:fixedpoint}, that is find $\vec x$ such that $\vec T(\vec x)=\vec 0$. An alternative formulation of \ac{aap} follows the definition of a fixed-point problem, that is find $\vec x$ such that $\vec x = \vec S(\vec x)$, for some given operator $\vec S$. Problems written in residual form allow one to compute the required vectors $\vec f^k$, $\vec g^k$ in the following order:
    \begin{equation}  \vec x^k \to \vec f^k = \vec T(\vec x^k) \to \vec g^k = \vec x^k - \omega \vec f^k \to \vec x^{k+1} = \vec g^k. \end{equation}
    If one considers an operator given in fixed-point form, then the second and third operations are executed in inverted order: 
    \begin{equation} \vec x^k \to \vec g^k = \vec S(\vec x^k)  \to \vec f^k = \vec x^k - \vec g^k \to \vec x^{k+1} = \vec x^k - \omega \vec f^k. \end{equation}
    Although our implementation in \texttt{Julia} supports both variants, the mathematical discussion presented henceforth relies only on the residual form, since most non-linear problems of interest in physics and engineering applications are described in residual form, thus justifying it to be the standard choice for most non-linear solver libraries. For readers interested in mathematical formulations of the \ac{aa} scheme for fixed-point iterations that commonly arise for solving non-linear \acp{pde}, we point to the \cite{pollock2019anderson} where \ac{aa} is presented and analyzed as an acceleration technique to speed up the Picard iterations for Navier-Stokes equations.
\end{remark}

\section{Sketching-based \ac{aap}}\label{section:aar sketching}
{To perform controlled perturbations of the Anderson least-squares problem whose admissibility is governed by the backward stability bounds derived in Section \ref{section:sketching analysis} of this work,} we consider two different levels of sketching in our algorithm, represented by restriction operators $\mat \Pi$ that project the solution onto a lower-dimensional space. To illustrate this approach, we focus on projections that extract a subset of indices, commonly referred to as \emph{masking} operators. 

The first level of restriction is physics based, and is primarily motivated by the multi-physics nature of certain problems such as the interaction between velocity and pressure in fluid dynamics. This approach is based on two key insights: (1) in many cases, one field can be effectively reduced in terms of another, which is rigorously justified in linear problems through the computation of the Schur complement matrix \cite{zhang2006schur}. This suggests that it is possible to solve the least squares problem using only one of the sub-physics without loss of information. (2) Differences in scale between the various physical components of a model can negatively affect the conditioning of the monolithic problem. We denote this subspace identified by the first level restriction as $\R^{\ell_1}$, with the corresponding restriction operator $\mat \Pi_1: \R^n \to \R^{\ell_1}$, where $\ell_1 \leq n$. This formulation allows us to determine the size of the reduced least squares matrix $\mat F_\Pi\coloneqq \mat{\Pi}_1 \mat F\in \R^{\ell_1\times n}$ at the allocation phase.

The second level of restriction is purely algebraic (not driven considerations on the physics) and is also dynamic (it can change across iterations). It is motivated by adaptive sketching techniques, such as the one proposed in \cite{lupo2022anderson}, based on the observation that precision in the least squares problem becomes less critical as the residual approaches zero. We define this second subspace as $\R^{\ell_2}$, with the corresponding restriction operator $\mat \Pi_2: \R^{\ell_1} \to \R^{\ell_2}$. Consequently, the overall restriction operator is given by $ \mat \Pi \coloneqq \mat \Pi_2 \circ \mat \Pi_1: \R^n \to \R^{\ell_2}$.

The algorithm that we propose thus extends Algorithm~\ref{alg:aar} with the following considerations: 
    \begin{enumerate}
        \item There is an \emph{allocation} phase where all required matrices and vectors are created.
        \item At each \ac{aap} iteration, there is a \emph{matrix update} phase where the newly computed residual is used to update the reduced matrix $\mat F_\Pi$ and the full matrix $\mat G$.
    \end{enumerate}
    The resulting algorithm for problem \eqref{sec2:fixedpoint} is shown in Algorithm~\ref{alg:new aar}, where the arrow ($\gets$) is used to represent \emph{in-place} BLAS operations (such as \texttt{axpy}). In this algorithm we have adopted the notation of writing all vectors without super indices, such as $\vec f$ instead of $\vec f^k$, to emphasize that we are reusing the same vector $\vec f$ instead of allocating a new one. This highlights the memory-reuse of our algorithm.

\begin{algorithm}
    \caption{Two-level sketching \ac{aap}}\label{alg:new aar}
    \begin{algorithmic}[1]
        \STATE{Memory allocation() \hfill (Algorithm~\ref{alg:allocation})}
        \STATE{Set $\vec x \gets \vec x^0$, $\vec f \gets \vec T(\vec x^0)$ and $\vec g \gets \vec x^0 + \omega \vec f^0$}
        \STATE{Set $\texttt{error}=1$, $\texttt{it}=0$}
        \WHILE{$\texttt{error} \geq \texttt{tol}$ and $\texttt{it} \leq \texttt{maxit}$}
            \STATE{Prepare $\Delta \vec f\gets \vec f$, $\Delta \vec g \gets \vec g$}
            \STATE{Compute residual $\vec f\gets \vec T(\vec x)$, $\vec g \gets \vec x - \omega \vec f$}
            \STATE{Update $\Delta \vec f\gets \vec f - \Delta \vec f$, $\Delta \vec g \gets \vec g - \Delta \vec g$}

            \IF{\texttt{adaptive}}
                \STATE{Update $\mat \Pi_2$ \hfill (Algorithm~\ref{alg:adaptive step})}
            \ENDIF

            \IF{$\mat\Pi_1\neq \mat I_n$}
                \STATE{Update $\vec f_\Pi \gets \mat \Pi_1\vec f$, $\Delta \vec f_\Pi \gets \mat\Pi_1 \Delta \vec f$ }
            \ENDIF

            \STATE{Delete first column from $\mat F_\Pi$ (if full) and shift columns}
            \STATE{Append $\Delta \vec g$ to $\mat G$ and $\Delta \vec f_\Pi$ to $\mat F_\Pi$}

            \IF{$k \neq 0 \mod p$}
                \STATE{Update $\vec x \gets \vec x - \omega\vec f$}
            \ELSE
            \STATE{Solve LS $\vec \alpha = \left[\mat\Pi_2 \mat F_\Pi\right]^\dag \mat\Pi_2\vec f_\Pi$ with $QR$ factorization, store $\mat R$}
                \STATE{Update $\vec x \gets \vec x - \omega\vec f + \mat G\vec \alpha$}
            \ENDIF
            \STATE{Update \texttt{error} and $\texttt{it}$}
        \ENDWHILE
        \RETURN $\vec x$
        \end{algorithmic} 
    \end{algorithm}

    The first step of the algorithm (line 1) involves memory allocation. We categorize the allocated objects based on their dimensions, as shown in Algorithm~\ref{alg:allocation}. We highlight that (i) we store only the matrix $\mat G$ and never $\mat X$, and (ii) if the sketching operator $\mat \Pi_1$ is used, the memory required for the matrix $\mat F$ can be reduced from $n$ rows to $\ell_1$ rows. Typically, the window size $m$ is small (5–10) and, in practice, it does not exceed 50–100. It is possible to create a copy of $\mat F_\Pi$ to enable a more memory-efficient implementation of the $QR$ factorization for the least squares (LS) problem. However, we have not observed any significant advantage in doing so, as it would require copying the entire (dense) matrix at each LS step. Therefore, we do not report the impact of this choice. We present a backward stability analysis of using restriction operators for the LS problem in Section~\ref{section:sketching analysis}. Additionally, the matrix $\mat R$ from the QR factorization is stored, as it is used during adaptivity to estimate the smallest singular value of the LS system.

    \begin{algorithm}[H]
        \caption{Memory allocation procedure}\label{alg:allocation}
        \begin{algorithmic}[1]
            \STATE{Allocate vectors $\vec x$, $\vec f$, $\vec g$, $\Delta \vec f$, $\Delta\vec g$ in $\R^n$}
            \IF{$\Pi_1 \neq \mat I_n$}
                \STATE{Allocate vectors $\vec r_\Pi$, $\Delta\vec r_\Pi $ in $\R^{\ell_1}$}
            \ENDIF
            \STATE{Allocate matrix $\mat G$ in $\R^{n\times m}$}
            \STATE{Allocate matrix $\mat F_\Pi$ in $\R^{\ell_1\times m}$}
            \STATE{If adaptive: allocate matrix $\mat R$ in $\R^{m\times m}$}
        \end{algorithmic} 
    \end{algorithm}

    When computing the residual (Algorithm~\ref{alg:new aar}, lines 5–7), we use an \emph{in-place} update of the increment vectors to save memory. The adaptivity step (Algorithm~\ref{alg:new aar}, line 9) is performed immediately afterward, as it only requires the residual vector (see Section~\ref{section:sketching analysis} for more details). The matrix update step (Algorithm~\ref{alg:new aar}, lines 14–15) follows these considerations: (i) the column deletion procedure is applied to the matrix $\mat F_\Pi$, highlighting once again the advantage of using a masking operator $\mat\Pi_1$, and (ii) the matrix $\mat G$ is updated in a circulant fashion, meaning that at iteration $k$, we update column $(k+1) \mod m$. This approach avoids shifting the entire matrix to the left, instead allowing direct access to the corresponding columns when updating the solution (Algorithm~\ref{alg:new aar}, line 20). 

\section{Projection-based sketching stability analysis}\label{section:sketching analysis}
In \ac{aap}, it is crucial to compute the Anderson mixing vector $\vec \alpha$ accurately. However, it has been observed that as iterations progress, the requirement for precision becomes increasingly less critical \cite{lupo2022anderson}. This observation can be rigorously justified through a backward stability analysis of the perturbed least squares problem, where the perturbations considered in practice arise from applying the restriction operators $\mat\Pi_1$ and $\mat\Pi_2$ to reduce the dimensionality of the LS problem, as described in Section~\ref{section:aar sketching}. 

In the following, we extend the proof provided in \cite{lupo2022anderson} to the non-linear case, under the assumption that the non-linear function $\vec T$ is Lipschitz continuous. That is, there exists a positive constant $L$ such that for all $\vec x_1, \vec x_2 \in \R^n$,
\begin{equation}\label{eq:hyp Lips}
    \| \vec T(\vec x_1) - \vec T(\vec x_2) \| \leq L \| \vec x_1 - \vec x_2 \|.
\end{equation}
Since restriction operators return vectors of reduced dimensionality, in the backward stability analysis we will use instead the induced projection operators, defined as $\mat P = \mat\Pi^T\mat\Pi$, $\mat P^1=\mat\Pi_1^T\mat\Pi_1$, and $\mat P^2 = \mat\Pi_2^T\mat\Pi_2$ respectively. These operators preserve the dimension of the vectors/matrices but zero-out the unused rows. 
For simplicity, we will drop the $\Pi$ subindex and thus use the matrix $\mat F^k$ instead of $\mat F^k_\Pi = \mat \Pi_1 \mat F^k$, so that on iteration $k$ the problem of interest yields the following solution: 
    \begin{equation}\label{eq:ls}
        \vec \alpha^k = \argmin_{\vec \alpha} \| \mat F^k \vec \alpha - \vec f^k \|^2. 
    \end{equation}
    Therefore, the optimal backward error to solve \eqref{eq:ls} is given by 
    \begin{equation}\label{eq:optimal backward error}
        \min_{\delta \mat F, \delta \vec f}\left\{\|\delta \mat F,\delta \vec f\|: (\delta \mat F, \delta \vec f) = \argmin_{(\delta \mat F, \delta \vec f)}\|(\mat F + \delta \mat F)\vec \alpha^k - (\vec f^k + \delta \vec f) \|^2 \right\}.
    \end{equation}
In particular, we assume that the error introduced on the right-hand side (RHS), represented by the vector $\delta \vec f$, can be controlled. Moreover, this error is used to bound the error introduced by $\delta\mat F$ on the left-hand side (LHS). Specifically, we employ projection operators $\mat P_k$ acting on each column, such that the perturbed matrix $\hat{\mat F}^k$ is given by
\begin{equation}\label{eq:perturbed matrix}
    \hat{\mat F}^k \coloneqq \mat F^k + \delta \mat F^k = \mat F^k - \left[ (\mat I - \mat P_1) \Delta \vec f^1, \hdots, (\mat I - \mat P_k)\Delta \vec f^k\right], 
\end{equation}
and similarly, we define the perturbed vector as
\begin{equation}\label{eq:perturbed vector}
    \hat{\vec f}^k \coloneqq \vec f^k + \delta \vec f^k = \vec f^k - (\mat I - \mat P_k)\vec f^k.
\end{equation}
    In what follows, we consider the $QR$ factorizations $\mat F^k = \mat Q^k\mat R^k$ and $\hat{\mat F}^k = \hat{\mat Q}^k\hat{\mat R}^k$ for the original and perturbed problems respectively.  We can now study the perturbation norm, defined as
    \begin{equation}\label{eq:perturbation norm}
        \delta^k = |\delta \mat F^k\hat{\vec\alpha}^k|_2 = \left| \left[(\mat I - \mat P_1) \Delta \vec f^1, \hdots, (\mat I - \mat P_k)\Delta \vec f^k\right] \hat{\vec \alpha}^k\right|_2 .
        \end{equation}
    Resorting to the Lipschitz continuity of the fixed-point operator $\vec T$, we obtain $|\Delta \vec f^k| =|\vec f^k - \vec f^{k-1}| \leq L | \Delta \vec x^k|$, and thus we can bound the perturbation as 
    \begin{equation}\label{eq:pert bound}
    \delta^k \leq \sum_i^k |\hat{\vec\alpha}^k|_2\|\mat I - \mat P_i \||\Delta \vec f^k|_2 \leq \sum_i^k L |\hat{\vec\alpha}^k|_2\|\mat I - \mat P_i \||\Delta \vec x^k|_2. 
    \end{equation} Using this estimate we can show the following result, which aids in devising adaptivity strategies for a judicious sketching in  Algorithm~\ref{alg:new aar}.

    \begin{theorem}\label{thm:adaptive} 
        Consider $\hat{\mat R}^k$ the R factor from the QR factorization of the perturbed matrix $\hat{\mat F}^k=\mat F^k + \delta \mat F^k$, $\epsilon$ a positive scalar representing the tolerable perturbation magnitude, $\vec f^k$ the $k$-th residual of \ac{aap}, and $(\eta_j)_j$ an arbitrary sequence such that $\sum_j^{\infty} \eta_j = C$ for some positive scalar $C$. We will further denote, for a given matrix $\mat Z$, $\sigma_\text{min}(\mat Z)$ as the smallest singular value of $\mat Z$. Consider the perturbations induced by the sketching operators as defined in \eqref{eq:perturbed matrix} and \eqref{eq:perturbed vector}, for which we assume the following hypotheses:
        \begin{enumerate}
            \item The perturbation $\delta \vec f^k$ of the k-th residual, defined by the relation $\hat{\vec f}^k = \vec f^k + \delta \vec f^k$, is such that 
                \begin{equation}
                    |\delta \vec f^k|_2 \leq \epsilon | \vec f^k|_2.
                \end{equation}
            \item The projection operators $(\mat I - \mat P_j)$ satisfy that
                \begin{equation}\label{eq:stability hyp}
                    \|\mat I  - \mat P_j \| \leq \eta_j\frac{\sigma_\text{min}(\hat{\mat R}^k)}{L|\vec f^k|_2|\Delta \vec x^j|_2(1+\epsilon)}\quad \forall j \in \{1,\hdots, k\}. 
                \end{equation}
        \end{enumerate}
        Under these assumptions, the perturbation norm $\delta^k$ defined in \eqref{eq:perturbation norm} satisfies that
            \begin{equation}
                \delta^k \leq C.
            \end{equation}
        \begin{proof}
            The proof closely resembles the one presented in \cite{lupo2022anderson}. We first bound the weights in terms of the perturbed vector. Using the relation $\left[\hat{\mat F}^k\right]^\dag = [\hat{\mat R}^k]^{-1}[\hat{\mat Q}^k]^T$, we obtain that
            \begin{multline}
                |\hat{\vec \alpha}^k|_2 = \left| [\hat{\mat R}^k]^{-1}[\hat{\mat Q}^k]^T\hat{\vec f}^k\right|_2 \\
                    \leq \left\|[\hat{\mat R}^k]^{-1}\right\|\left|\hat{\vec f}^k\right|_2 \leq \frac{1}{\sigma_\text{min}(\hat{\mat R}^k)}\left(\left|\vec f^k\right|_2 + \left|\delta \vec f^k\right|_2\right) \leq \frac{1+\epsilon}{\sigma_\text{min}(\hat{\mat R}^k)}\left|\vec f^k\right|_2, 
            \end{multline}
            where we used that $\|\mat{Z}^{-1}\| = \frac{1}{\sigma_\text{min}(\mat Z)}$. Putting this in \eqref{eq:pert bound} we obtain  \begin{equation}\label{eq:pert bound 2}
                    \delta^k \leq \sum_i^k \frac{L(1+\epsilon)}{\sigma_\text{min}(\hat{\mat R}^k)}\left|\vec f^k\right|_2\left\|\mat I - \mat  P_i\right\| \left|\Delta x^i\right|_2.
                \end{equation}
            From this equation we can then conclude that, if for every $j\leq k$ \eqref{eq:stability hyp} holds, then  we have 
                \begin{equation} 
                    \delta^k \leq \sum_j^k \eta_j \leq C,
                    \end{equation}
            which establishes the result.
        \end{proof}
    \end{theorem}
    Using this theoretical result, we can provide simple but still effective strategies to adaptively project the LS problem onto a subspace defined by reducing the \acp{dof} retained from the original LS problem.  Designing an adaptive strategy needs to consider, on one hand, Theorem~\ref{thm:adaptive} in order to ensure that the approximated calculations to solve the LS problem at each mixing step do not affect the convergence of \ac{aap} to the desired accuracy. On the other hand, it needs to consider the slowdown associated to non-sequential memory access due to a highly frequent cache miss. We have studied what is the maximum percentage of sketching that still yields faster matrix-vector products and faster QR factorizations in Section~\ref{sec:sketching-benchmarks} with respect to a dense matrix with sequential memory access.
    
    \subsection{Ensuring LHS stability}
    We start by observing that, as $\mat I - \mat P_k$ are projection operators, they have unit norm. Therefore, the second hypothesis on Theorem~\ref{thm:adaptive} becomes 
        \begin{equation}\label{eq:restriction epsilon intermediate}
            \|\mat I  - \mat P_j \| = 1 \leq \eta_j\frac{\sigma_\text{min}(\hat{\mat{R}^k})}{L|\vec f^k||\Delta \vec x^j|(1+\epsilon)} \quad \forall j \in \{1,\hdots, k\},
        \end{equation}
        which can be recast as a bound on $\epsilon$ as it must hold that
        \begin{equation}\label{eq:epsilon bound}
            \epsilon \leq \frac{\eta_j\sigma_\text{min}(\hat{\mat{R}^k})}{L|\vec f^k||\Delta \vec x^j|} - 1=: \epsilon_\text{LHS}.
        \end{equation}

        Since $\sigma_\text{min}(\hat{\mat R}^k)$ and $L$ are impractical to compute, we adopt the following computationally inexpensive approximations:
        \begin{itemize}
            \item Estimate $\sigma_\text{min}$ by performing a few (e.g., between 1 and 5) iterations of the inverse power method to the matrix $[\hat{\mat R}^k]^T\hat{\mat R}^k$ \cite{kelley1995iterative}.
            \item Estimate $L$ by keeping track of the increment norms so that
                \begin{equation}\label{eq:L estimate}
                    L \approx L_k = \max\left(L_{k-1},  \frac{|\Delta \vec f^k|_2}{|\Delta \vec x^k|_2}\right).
                \end{equation}
        \end{itemize}

    \subsubsection{Ensuring RHS stability}
        To enforce this constraint, we will use that the perturbation vector $\delta \vec f^k = (\mat I - \mat P^k)\vec f^k$ is given by the entries that are removed from the vector, and ensure that $|\delta \vec f^k| < \epsilon |\vec f^k|$. An effective way to use this is to compute some indices, and with that simply define the quantity
        \begin{equation}\label{eq:erhs definition}
            \epsilon_\text{RHS} \coloneqq \frac{|\delta\vec f^k|_2}{|\vec f^k|_2}. 
        \end{equation}
        To establish a connection between $\epsilon_\text{LHS}$ and $\epsilon_\text{RHS}$, we impose that $\epsilon = \epsilon_\text{RHS}$. From Equation~\ref{eq:epsilon bound} we have that $0<\epsilon_\text{RHS} \leq \epsilon_\text{LHS}$.

    \subsection{Choosing the index set}\label{sec:sketching-benchmarks}
\begin{figure}
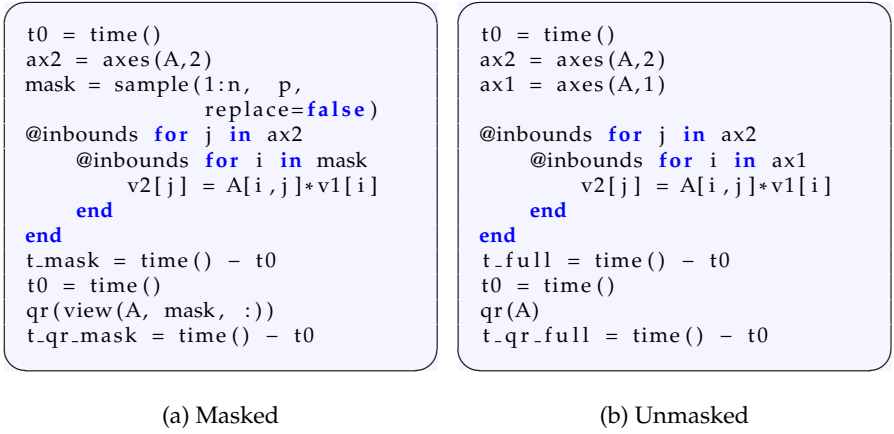

        \centering
        \begin{subfigure}[t]{0.4\textwidth}
            \centering
            \begin{lstlisting}[language=Julia]
t0 = time()
ax2 = axes(A,2)
mask = sample(1:n,  p, 
              replace=false)
@inbounds for j in ax2
    @inbounds for i in mask
        v2[j] = A[i,j]*v1[i]
    end
end
t_mask = time() - t0
t0 = time()
qr(view(A, mask, :))
t_qr_mask = time() - t0
            \end{lstlisting}
            \caption{Masked}
        \end{subfigure}
        \begin{subfigure}[t]{0.05\textwidth}
        \phantom{asdf}
        \end{subfigure}
        \begin{subfigure}[t]{0.4\textwidth}
            \centering
            \begin{lstlisting}[language=Julia]
t0 = time()
ax2 = axes(A,2)
ax1 = axes(A,1)

@inbounds for j in ax2
    @inbounds for i in ax1
        v2[j] = A[i,j]*v1[i]
    end
end
t_full = time() - t0
t0 = time()
qr(A)
t_qr_full = time() - t0
            \end{lstlisting}

            \caption{Unmasked}
        \end{subfigure}
        \caption{Matrix-vector product and QR computation benchmarks (a) with and (b) without masking.}
        \label{fig:view code}
    \end{figure}

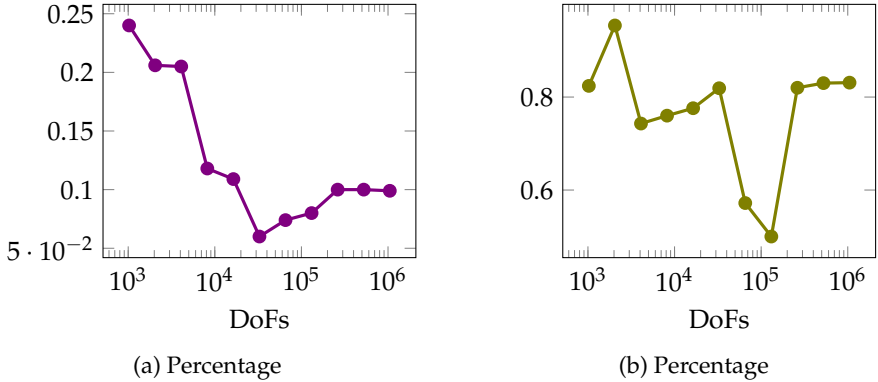
\begin{figure}
        \centering
        \begin{subfigure}{0.49\textwidth}
            \centering
            \begin{tikzpicture}
                \begin{semilogxaxis}
                    [xlabel=\acp{dof},
                    width=0.9\textwidth
                    ]
                    \addplot[mark=*, line width=1.2, blue!50!red] coordinates {
                                          (1024   , 0.240) 
                                          (2048   , 0.206)
                                          (4096   , 0.205)
                                          (8192   , 0.118)
                                          (16384  , 0.109)
                                          (32768  , 0.060)
                                          (65536  , 0.074)
                                          (131072 , 0.080)
                                          (262144 , 0.100)
                                          (524288 , 0.100)
                                          (1048576, 0.099)
                                          };
                \end{semilogxaxis}
            \end{tikzpicture}
            \caption{Percentage}
        \end{subfigure}
        \begin{subfigure}{0.49\textwidth}
            \centering
            \begin{tikzpicture}
                \begin{semilogxaxis}
                    [xlabel=\acp{dof},
                    width=0.9\textwidth
                    ]
                    \addplot[mark=*, line width=1.2, green!50!red] coordinates {
                                          (1024   , 0.824) 
                                          (2048   , 0.954)
                                          (4096   , 0.743)
                                          (8192   , 0.760)
                                          (16384  , 0.776)
                                          (32768  , 0.819)
                                          (65536  , 0.572)
                                          (131072 , 0.500)
                                          (262144 , 0.820)
                                          (524288 , 0.830)
                                          (1048576, 0.831)
                                          };
                \end{semilogxaxis}
            \end{tikzpicture}
            
            \caption{Percentage}
        \end{subfigure}
        \caption{Maximum allowed percentage of random indices required to have faster (a) matrix vector product and (b) QR factorization, for increasing vector sizes.}\label{fig:mask percentage}
    \end{figure}

    To determine the rows to remove in the tall-and-skinny matrix $\mat F^k_\Pi$ and its corresponding right hand side $\vec f^k$, we consider that, in practical terms, accessing random indices can be slower than accessing the full vector due to frequent cache misses in the random access memory of the computer. For this reason, we estimated the maximum number of random indices that can be used while still ensuring a faster matrix-vector product and QR factorization. We present the benchmarking code to compare the computational time for matrix-vector product and QR factorization with and without randomized sketching in Figure~\ref{fig:view code}. We performed the benchmarking using a random matrix with 50 columns and a number of rows ranging from $10^3$ up to $10^6$, and show the maximum allowed percentages in Figure~\ref{fig:mask percentage}. To make these estimates precise, we ran every test and averaged the results until the relative change in the computed average was below $10^{-5}$. 

From this benchmarking test, we conclude that: our sketching operators should retain no more than 10\% of the original rows to achieve a faster matrix-vector product and no more than 80\% of the original rows to ensure a faster QR factorization. Following the approach in \cite{lupo2022anderson}, we consider two sketching strategies for a given percentage $S$:
\begin{itemize}
    \item \textbf{Subselection}: We retain the $S\%$ of entries corresponding to the largest current residual $\vec f^k$.
    \item \textbf{Randomized}: We retain a randomly selected $S\%$ of the entries.
\end{itemize}
The subselection strategy follows the intuition that the portions of the residual with larger residuals are where computational efforts should be concentrated to effectively reduce the error.

    \subsection{Overall strategy}\label{section:adaptive step}
    Under all previous considerations, the final adaptivity strategy requires the following steps: 
        \begin{enumerate}
            \item Fix a sketching percentage $S$, 
            \item estimate  $\epsilon_\text{LHS}$ and assert it is positive, 
            \item compute the index set defining the mask, 
            \item use the mask to compute $\epsilon_\text{RHS}$, and 
            \item if $0\leq\epsilon_\text{RHS}\leq \epsilon_\text{LHS}$, use mask. Otherwise, do not use $\mat\Pi_2$.
        \end{enumerate} 
        We show the detailed steps in Algorithm~\ref{alg:adaptive step}. We have added the problem dimension $N$ to the $\epsilon_\text{LHS}$ estimator in order to consider vector norms that are less prone to dimensionality problems with $|\cdot |^2 = \frac 1 N \sum_i^N (\cdot)^2$. Without this, adaptivity is seldom activated for larger problems. A more general strategy to avoid this kind of dimensionality problems would be to frame the method in terms of the infinite dimensional spaces where the \acp{pde} are formulated and compute the corresponding infinite dimensional norms. This is beyond the scope of this work and part of what we plan to address in future work. 

        We consider two choices for the subselection parameter $\eta_j$: either $\eta_j = j^{-1.1}$ (termed \emph{power} adaptivity) or $\eta_j = 1$ (termed \emph{constant} adaptivity). Combining these variants yields five adaptivity strategies used in our numerical experiments:
\begin{itemize}
  \item \textbf{No adaptivity}: $\mat \Pi_{2} = \mat I$.
  \item \textbf{Subselect power}: Residual subselection with $\eta_j = j^{-1.1}$.
  \item \textbf{Subselect constant}: Residual subselection with $\eta_j = 1$.
  \item \textbf{Randomized power}: Random index selection with $\eta_j = j^{-1.1}$.
  \item \textbf{Randomized constant}: Random index selection with $\eta_j = 1$.
\end{itemize}

    Our analysis requires that series given by $(\eta_j)_j$ is convergent, and thus the disctinction between the \emph{constant} and \emph{power} algorithms is meaningful as only the latter is strictly within our theoretical setting. We have chosen somewhat arbitrarily the exponent 1.1 as it yields a convergent series with a slow decrease. Larger exponents are naturally possible but they could speed-up the decrease of the term $\eta_j \sigma_\text{min}(\hat{\mat R}^k)/L|\vec f^k||\Delta x^j|$ and thus artificially prevent adaptivity. We present a sensitivity study of adaptivity on the exponent $\ell$ of $\eta^{-\ell}$ in Appendix \ref{appendix:power}, where we show for a specific problem that it is not significative for the random masks, and negligible for subselection.

    \begin{algorithm}
        \caption{Adaptive step for computing $\mat\Pi_2$}\label{alg:adaptive step}
        \begin{algorithmic}[1]
            \STATE{Compute $\sigma \approx \sigma_\text{min}(\hat{\mat R}^k)$ using few (1-5) inverse iterations.}
            \STATE{Compute $\epsilon_\text{LHS} =  \max_j \frac{N\eta_j\sigma}{L|\vec f^k||\Delta \vec x^j|} - 1$.}
            \IF{$\epsilon_\text{LHS} < 0$}
                \STATE{\textbf{break}}
            \ELSE
                \STATE{Fix number of rows $\ell_2$ using $S$.}
                \STATE{Compute index set $\mathcal I^j$ (subselection or randomized).}
                \STATE{Compute $\epsilon_\text{RHS} = |(\mat I - \mat\Pi)\Delta\vec f^j|/|\vec f^j|$.}
                \IF{$0<\epsilon_\text{RHS} \leq \epsilon_\text{LHS}$}
                    \STATE{$[\mat P^j_2]_{ik} = \begin{cases} 1 & \text{if $i=k \in \mathcal I^j$} \\ 0\end{cases}$}
                \ELSE
                    \STATE{$\mat P^j_2=\mat I$}
                \ENDIF
            \ENDIF
        \end{algorithmic} 
    \end{algorithm}

\section{Numerical tests}\label{section:tests}
In this section, we apply Algorithm~\ref{alg:new aar} to solve several different \acp{pde}. These numerical tests aim to show that the proposed methodologies remain robust across varying problem sizes and can yield improvements in overall solution time. All tests were implemented in Julia \cite{Bezanson_Julia_A_fresh_2017} using the Ferrite \cite{CarlssonFerritejl} and Gridap \cite{badia2020gridap} finite element libraries, and were executed on a single core of an AMD Ryzen Threadripper 3990X. To ensure accurate timing measurements, tests were run one at a time, with all non‑essential Linux kernel services disabled during execution and the compiler optimization flags \texttt{-O3 --check‑bounds=no} enabled. The problems we studied are:
\begin{enumerate}
    \item Stokes: a symmetric but indefinite problem coupling velocity and pressure.
    \item $p$-Laplacian: a nonlinear problem that becomes challenging when $p$ deviates significantly from 2.
    \item Bidomain equations: a degenerate parabolic system modeling the propagation of electrical potentials in cardiac tissue.
    \item Stationary Navier–Stokes equations for the lid-driven cavity test.
\end{enumerate}

For each \ac{pde} problem, we performed the following tests:
    \begin{itemize}
        \item \textbf{Physics-based mask test.} We propose possible mask operators $\mat \Pi_1$ according to the sub-blocks of the problem and use them with \ac{aap}. We report the iteration number and CPU time in seconds for increasing problem size.
        \item \textbf{Adaptivity test.} We study all five considered adaptivity strategies (no adaption, subselection with power sequence, subselection with constant sequence, randomized with power sequence, and randomized with constant sequence). In each case, we report the iteration count and the CPU time in seconds for increasing problem size.
        \item \textbf{Best overall performance.} We test all combinations of alternation parameters $p$ in $\{1,2,3,4\}$, physics-based masks, and adaptivity strategies. For each problem size considered in the previous tests, we show the method with the best performance as measured by the CPU time. 
    \end{itemize}
    Unless otherwise stated, we used a window size of $m=10$ and no alternation $p=1$ as a base method for comparison, and a sketching percentage of $30\%$ for all adaptive strategies.

\subsection{Stokes problem}
The Stokes problem is formulated as follows: Given $f:\Omega \to \R^3$, find a velocity $\vec u:\Omega \to \R^3$ and pressure $\pressureS:\Omega\to \R$ such that 
    \begin{equation}
        \begin{aligned}
            -\Delta \vec u + \grad \pressureS &= \vec f &&\text{ in $\Omega$},\\
            \dive \vec u &= 0 &&\text{ in $\Omega$},\\
            \vec u &= \vec 0 &&\text{on $\Gamma_D$}, \\
        \left(\grad \vec u -\pressureS\ten I\right)\vec n &= \vec 0 &&\text{on $\Gamma_N$}.
        \end{aligned}
    \end{equation}
    We considered $\Omega=(0,1)^3$, a constant volume force $\vec f(x)=(1,1,1)$, homogeneous Dirichlet conditions on $\Gamma_D=\{\vec x=(x,y,z)\in \partial\Omega:  y\in\{0,1\}, z\in \{0,1\}\}$, and homogeneous Neumann conditions on $\Gamma_N=\{\vec x=(x,y,z)\in \partial\Omega: x\in\{0,1\}$.  This problem is discretized with inf-sup stable Taylor-Hood elements $\mathbb P_2\times \mathbb P_1$ \cite{elman2014finite}, which gives rise to the block linear system 
        \begin{equation}
            \begin{bmatrix} \mat K & \mat B^T \\ \mat B & \mat 0 \end{bmatrix}\begin{bmatrix} \vec u\\ \pressureS \end{bmatrix} = \begin{bmatrix}\vec F \\ \vec 0 \end{bmatrix}.
        \end{equation}
        We use a preconditioned Richardson iteration
        \begin{equation}
            \begin{bmatrix}\vec u^{k+1} \\ \pressureS^{k+1}\end{bmatrix} = \begin{bmatrix}\vec u^{k} \\ \pressureS^{k}\end{bmatrix} - \mat P^{-1}\left( \begin{bmatrix} \mat K & \mat B^T \\ \mat B & \mat 0 \end{bmatrix}\begin{bmatrix} \vec u^k\\ \pressureS^k \end{bmatrix} - \begin{bmatrix}\vec F \\ \vec 0 \end{bmatrix}
            \right),
        \end{equation}
        and accelerate it with our implementation of \ac{aap}. We use a block diagonal preconditioner $\mat P$ with the stiffness matrix on the velocity block and a mass matrix in the pressure block. This is known to be a scalable preconditioner for this problem \cite{elman2014finite}, and we approximate the inverses with an Algebraic Multigrid for both blocks as provided by HYPRE \cite{falgout2002hypre}. We study the effect of physics-based masks and algebraic masks obtained using adaptivity in terms of number of iterations and time-to-solution, and we present also the performance of the approach identified as best performing among all the combinations of masking and adaptive selection that we tested. A relative tolerance of $10^{-6}$ is used to determine convergence. 

    \subsubsection{Stokes: Physics-based mask test}
    The physics-based restriction operator we consider is the one induced by the physics of the problem, meaning either
        \begin{equation}\label{eq:up-masks}
            \mat \Pi_1\left(\begin{bmatrix} \vec u \\ \pressureS \end{bmatrix}\right) = \vec u \qquad \text{ or } \qquad  \mat \Pi_1\left(\begin{bmatrix} \vec u \\ \pressureS \end{bmatrix}\right) = \pressureS,
        \end{equation}
        which we refer to as \emph{velocity} and \emph{pressure} masks, respectively. We show in Table~\ref{table:stokes-mask} the total iterations and solution time in seconds required to converge with $\mat\Pi_1$ being an identity (no mask), a velocity mask, and a pressure mask. 

We note that the velocity mask yields more iterations and higher time-to-solution than the unmasked case, whereas the pressure mask converges in fewer iterations and in less time than the unmasked case. More specifically, the use of the pressure mask results in a 24\% reduction in the number of iterations and a 33\% reduction in the solution time in the largest problem when compared to the unmasked version. The conclusion we draw from these results is that the pressure seems to capture the complexity of the problem more effectively than the velocity, which seems to be a reasonable conclusion given that the Stokes equation can be reformulated as a pressure-only problem by means of the Schur complement.

    \begin{table}[ht!]
      \centering
      \footnotesize
      \begin{tabular}{r | r r | r r | r r}
          \toprule \multirow{2}{0.8cm}{\acp{dof}} & \multicolumn{2}{c|}{No mask} & \multicolumn{2}{c|}{Velocity mask} & \multicolumn{2}{c}{Pressure mask} \\ 
                  & Iters & Time   &  Iters & Time   & Iters & Time   \\ \midrule
          1148    & 138   & 0.10   & 143    & 0.09   & 118   & {\bf 0.06}   \\
          12204   & 117   & {\bf 0.95}   & 158    & 1.31   & 149   & 1.06   \\
          100052  & 119   & 11.65  & 139    & 13.50  & 102   & {\bf 8.92}   \\
          809892  & 113   & 102.61 & 160    & 140.57 & 96    & {\bf 77.32}  \\ 
          2743924 & 119   & 385.12 & 128    & 409.25 & 91    & {\bf 259.46} \\ \bottomrule
      \end{tabular}
        \caption{Stokes problem, physics-based mask test: Number of iterations and CPU time in seconds required for convergence using no mask, a velocity mask, and a pressure mask for $\mat \Pi_1$.The lowest CPU time for each row has been highlighted with bold font.}\label{table:stokes-mask}
    \end{table}

    \subsubsection{Stokes: Adaptivity test}
    Table~\ref{fig:stokes-adapt} shows the total iterations and CPU time in seconds required to solve the Stokes problem using all five adaptivity strategies described. As in the physics-based mask test, all strategies considered behave robustly with respect to the problem size in terms of the number of iterations. 

    Subselection strategies (power and constant) are stable as the number of iterations closely resembles that of the non adapted method. Randomized strategies present worse performance than the non adapted strategy for the smaller problems, as it incurs in more iterations to achieve convergence. 

    The Randomized Constant method requires more iterations than the other adaptive strategies, which is consistent with the theory as it uses the non-convergent series $\eta_j=1$ in \eqref{eq:epsilon bound}. Counterintuitively, this is not true for the larger problems, as in the second to last row the Randomized Constant strategy is comparable to subselection. We observe that all adaptivity strategies are faster than the non adapted one for problems with the largest size. In addition, Randomized Power is the most efficient strategy for the largest problem, which yields the lowest iterations and fastest execution with a reduction of CPU time by 17\% and a reduction of the number of iterations by 23\% compared to the non-adaptive approach. 

    \begin{table}[ht!]
        \centering
        \footnotesize
      \begin{tabular}{r | r r | r r | r r | r r | r r}
          \toprule \multirow{2}{0.8cm}{\acp{dof}} & \multicolumn{2}{c|}{No adapt} & \multicolumn{2}{c|}{Sub pow} & \multicolumn{2}{c|}{Sub const} & \multicolumn{2}{c|}{Rand pow} & \multicolumn{2}{c}{Rand const} \\ 
                  & Iters &  Time   &  Iters  & Time   &   Iters & Time   & Iters & Time   & Iters & Time  \\ \midrule
          1148    & 138   & {\bf 0.10}    & 140     & 0.13   & 131     & 0.12   & 130   & 0.12   & 243   & 0.22  \\
          12204   & 117   & 0.95    & 106     & {\bf 0.88}   & 108     & 0.92   & 174   & 1.44   & 188   & 1.55  \\
          100052  & 119   & 11.65   & 118     & {\bf 11.06}  & 118     & 11.07  & 137   & 12.68  & 143   & 13.58 \\
          809892  & 113   & 102.61  & 112     & {\bf 96.24}  & 112     & 96.93  & 126   & 106.47 & 116   & 97.85 \\
          2743924 & 119   & 385.12  & 104     & 315.79 & 104     & 318.87 & 99    & {\bf 297.57} & 110   & 329.95\\ \bottomrule
      \end{tabular}

      \caption{Stokes problem, adaptivity test: Number of iterations and CPU time in seconds required to solve the Stokes problem using no adaptivity, Subselection power,  Subselection constant,  Randomized power, and  Randomized constant for increasing problem size. The lowest CPU time for each row has been highlighted with bold font.}\label{fig:stokes-adapt}
      \end{table}

      \subsection{Stokes: Best overall performance test}
      Here we tested all combinations of the considered physics-based masks, adaptivity strategies, and alternation parameters ($p$ in $\{1,2,3,4\}$). For each problem size, we report in Table~\ref{table:stokes-best} the iterations and solution time of \ac{aap}(10,1) (equivalently \ac{aa}(10)) and the best performing method (in terms of CPU time). We note that there is no scenario in which standard \ac{aa} performs similarly to the sketching variants: iteration counts are significantly reduced in all masked cases and the same holds for CPU time, which is reduced by up to a 54\% in the largest problem. Even though we could not identify an adaptivity strategy that clearly outperforms consistently all the others, we noticed that every adaptive strategy performs better when combined with alternation, with a consistent value of $p=4$ being better for all cases considered.

    \begin{table}[ht!]
      \centering
      \footnotesize
      \begin{tabular}{r | rr | r r | l l c }
          \toprule \multirow{2}{0.8cm}{\acp{dof}} & \multicolumn{2}{c|}{No mask} & \multicolumn{5}{c}{Best}  \\ 
                                             & Iters        &  Time   &  Iters &    Time      &  Mask & Adaptivity & Alt ($p$)  \\ \midrule
          1148    & 138   & 0.10    & 84     & {\bf 0.04}   & Pressure mask& No adaptivity& $4$       \\
          12204   & 117   & 0.95    & 72     & {\bf 0.55}   & No mask& Subselect power& $4$    \\
          100052  & 119   & 11.65   & 68     & {\bf 6.13}   & No mask& Subselect constant& $4$ \\
          809892  & 113   & 102.61  & 64     & {\bf 53.25}  & No mask& Random power& $4$   \\ 
          2743924 & 119   & 385.12  & 60     & {\bf 175.30} & No mask& Random constant& $4$   \\ \bottomrule
      \end{tabular}
      \caption{Stokes problem, best overall test: Number of iterations and CPU time for the reference \ac{aa} method and the best performing one, characterized by its mask, adaptivity strategy, and alternation parameter $p$. The lowest CPU time for each row has been highlighted with bold font.}\label{table:stokes-best}
    \end{table}

\subsection{$p$-Laplace}
The $p$-Laplacian problem is given by: Given $\pL\geq 1$ and $f:\Omega\to\R$, find $u:\Omega \to \R$ such that
    \begin{equation}
        \begin{aligned}
            F(u) \coloneqq -\dive |\grad u|^{\pL-2}\grad u - f &= 0 &&\text{in $\Omega$}, \\
            u &= 0 &&\text{in $\partial \Omega$},
        \end{aligned}
    \end{equation}
    where $\Omega$ is the unit cube, $f(x)=1$, and $\pL=1.5$. We used $\pL$ instead of $p$ to avoid confusion with the alternation parameter $p$ in \ac{aap}. We chose $\pL<2$ as it is a challenging scenario for which there are fewer efficient strategies in the literature \cite{loisel2020efficient}. To solve the problem, we discretized it with first order Lagrangian finite elements and then used a quasi-Newton scheme to solve the nonlinear problem. The solver is given as follows: Consider a harmonic initial guess $u^0$ and then iteratively solve a Laplace-stabilized iteration:
    \begin{equation} -\beta \Delta \delta u^k = - F(u^{k-1}), \end{equation}
    where we have found the value $\beta=10$ to behave robustly. The increment $\delta u^k$ is used as a residual, i.e. 
    \begin{equation}\label{eq:plaplace-fp} \vec T(u^k) = \delta u^k = \frac 1 \beta  \Delta^{-1}\left(F(u^{k-1})\right), \end{equation}
    that allows for acceleration using \ac{aap}. Instead of inverting the Laplacian exactly in \eqref{eq:plaplace-fp}, we consider the action of one iteration of a smoothed aggregation AMG preconditioner, as provided by the PCGAMG preconditioner in PETSc \cite{petsc-user-ref}. For this problem, we test the adaptivity strategies and then look for the overall best performing method. There are no sub-physics in this problem, so we simply set $\mat\Pi_1=\mat I$, and a relative tolerance of $10^{-6}$ was used to establish convergence.

    \subsubsection{$p$-Laplace: Adaptivity test}
    We show the iteration count and CPU time in seconds required to solve this problem with all the proposed adaptivity strategies in Table~\ref{table:plaplace-adapt}. As in the previous tests, all formulations are robust with respect to the increase in \acp{dof} as seen by the iteration counts. In addition, iteration counts are similar to the non adapted case, which is consistent with the backwards stability analysis that inspired the formulation of these methods. 

    The subselection strategies perform very similarly to the non-adapted case in terms of CPU time, and yield a mild increase in the number of iterations which results in them being mildly slower in general. For the three largest problems, randomized adaptivity yields better performance (last two columns), with Randomized Constant yielding an iteration count and CPU time reduction of 13\% and 19\% respectively.

    \begin{table}[ht!]
        \centering
        \footnotesize
      \begin{tabular}{r | r r | r r | r r | r r | r r}
          \toprule \multirow{2}{0.8cm}{\acp{dof}} & \multicolumn{2}{c|}{No adapt} & \multicolumn{2}{c|}{Sub pow} & \multicolumn{2}{c|}{Sub const} & \multicolumn{2}{c|}{Rand pow} & \multicolumn{2}{c}{Rand const} \\ 
                   & Iters &  Time  & Iters  & Time   & Iters  & Time   & Iters & Time   & Iters & Time  \\ \midrule
          343      & 39    & {\bf 0.04}   & 44     & 0.06   & 49     & 0.07   & 41    & 0.06   & 40    & 0.05 \\
          3375     & 36    & {\bf 0.17}   & 41     & 0.20   & 39     & 0.19   & 40    & 0.20   & 39    & 0.19 \\
          29791    & 48    & 1.91   & 44     & {\bf 1.60}   & 53     & 1.91   & 49    & 1.76   & 49    & 1.75 \\
          250047   & 50    & 15.11  & 57     & 15.62  & 56     & 15.93  & 52    & 14.14  & 52    & {\bf 14.05} \\ 
          2048383  & 53    & 130.05 & 57     & 132.22 & 57     & 131.26 & 50    & 114.12 & 46    & {\bf 105.38} \\ \bottomrule
      \end{tabular}

      \caption{$p$-Laplace adaptivity test: Iteration count and CPU time in seconds for all adaptivity strategies proposed (none, Subselect power, Subselect constant, Randomized power, and Randomized constant) and for varying problem size. The lowest CPU time for each row has been highlighted with bold font.}\label{table:plaplace-adapt}
      \end{table}

    \subsection{$p$-Laplace: Best overall performance test}
    To obtain the best performing method in terms of CPU time, we considered all combinations of adaptivity strategies and alternation parameter $p$ in $\{1,2,3,4\}$. We display the results in Table~\ref{table:plaplace-best}, where for each method we show the iteration count and the CPU time in seconds. In this problem, we note that the best performing methods use either alternation, adaptivity, or both, with no particular trend. In all cases, the use of these strategies yields a reduction in the number of iterations and CPU time, where the best strategy in the largest problem is to simply use a Randomized Constant adaptivity strategy, as displayed in the previous test.

    \begin{table}[ht!]
      \centering
      \footnotesize
      \begin{tabular}{r | rr | r r | l c}
          \toprule \multirow{2}{0.8cm}{\acp{dof}} & \multicolumn{2}{c|}{No mask} & \multicolumn{4}{c}{Best}  \\ 
                  & Iters  &  Time  & Iters  &    Time     & Adaptivity & Alt ($p$)  \\ \midrule
          343     & 39     & 0.04   & 22     & {\bf 0.02}  & No adaptivity & $4$   \\
          3375    & 36     & 0.17   & 26     & {\bf 0.12}  & Random power & $3$ \\
          29791   & 48     & 1.91   & 40     & {\bf 1.42}  & Subselect Constant & $2$ \\
          250047  & 50     & 15.11  & 36     & {\bf 10.00} & No adaptivity  & $2$ \\
          2048383 & 53     & 130.05 & 46     & {\bf 105.38}& Random Constant & $1$ \\ \bottomrule
      \end{tabular}
      \caption{$p$-Laplace, best overall test. Number of iterations and CPU time for the reference \ac{aa} and the best performing method. The best methods are characterized by their adaptivity strategy and alternation parameter $p$. The lowest CPU time for each row has been highlighted with bold font.}\label{table:plaplace-best}
    \end{table}

\subsection{Bidomain equations}
The Bidomain model is given by the following system: Given an ionic model $I_\text{ion}:\Omega\to \R$ and applied current $I_\text{app}:\Omega \to \R$, find extracellular and intracellular potentials $u_e,u_i:\Omega \to \R$  such that
    \begin{equation}\label{eq:bidomain}
        \begin{aligned}
            F_e(u_e, u_i) \coloneqq \parder{(u_e - u_i)}{t} - \dive D_e\grad u_e + I_\text{ion}(u_e-u_i) - I_\text{app} &= 0  &&\text{in $\Omega$}, \\
            F_i(u_e, u_i) \coloneqq \parder{(u_i - u_e)}{t} - \dive D_i\grad u_i - I_\text{ion}(u_e-u_i) + I_\text{app} &= 0 &&\text{in $\Omega$}, \\
            \grad u_{i,e} \cdot \vec n &= 0 &&\text{on $\partial\Omega$},
        \end{aligned}
    \end{equation}
    where $D_{e,i}$ are symmetric and positive definite conductivity tensors. This is a degenerate parabolic problem that is commonly used to describe the propagation of an electric potential in soft tissue, such as the one that generates a heartbeat in cardiac tissue. We consider the same parameter setting as in \cite{barnafi2024robust}, with a Fitz-Hugh Nagumo ionic current. We discretize it with an implicit Euler scheme in time, and with first order finite elements in space. All the numbers we report are related to the solution of the first time step.  It has been observed in \cite{barnafi2024robust} that \ac{aa} yields a robust and scalable solver for this model.

To solve it, we consider the fixed point map given by the following Picard iteration:
    \begin{equation}
        \begin{bmatrix} u_e^{k+1} \\ u_i^{k+1} \end{bmatrix} = \begin{bmatrix} u_e^{k} \\ u_i^{k} \end{bmatrix} - \begin{bmatrix} F_e(u_e^k, u_i^k) \\ F_i(u_e^k, u_i^k) \end{bmatrix}
\end{equation}
and accelerate it with \ac{aap} with a window size of $m=50$. We have observed that smaller window sizes do not yield a convergent scheme. A relative tolerance of $10^{-4}$ was used to establish convergence.

    \subsubsection{Bidomain: Physics-based mask test}
    The structure of the Bidomain equations yields the following projection operators:
    \begin{equation} \mat \Pi_1\left(\begin{bmatrix} u_e \\ u_i \end{bmatrix}\right) = u_e, \qquad \text{ and } \qquad  \mat \Pi_1\left(\begin{bmatrix} u_e \\ u_i \end{bmatrix}\right) = u_i, \end{equation}
    which we refer to as \emph{extracellular} and \emph{intracellular} masks, respectively. We solved the problem using an identity, extracellular, and intracellular masks, and report the iteration count and CPU time in seconds obtained in Table~\ref{table:bidomain-mask}. In this test, it is important to note that the number of \acp{dof} that the extracellular and intracellular potentials have is the same, which makes their difference in performance particularly unexpected.

We start by observing that none of the physics-based masks converge for problems of small size. For larger problems (from 500k \acp{dof} onwards), the extracellular mask behaves similarly to the unmasked case, whereas the intracellular mask presents a performance deterioration that yields more than three times the number of iterations, and divergence in the largest problem. 
For a sufficiently large problem, both the unmasked and extracellular mask methods yield a robust behavior in terms of total iteration count, and this allows for the extracellular mask to have a similar (but mildly faster) performance than the unmasked case. CPU time savings in the largest problem considered are not significant (1\%). 

However, since both equations are similar and both fields have the same number of \acp{dof}, one would expect similar performance when adopting either one of the two variables for masking. Therefore, the performance deterioration when using a intracellular mask warrants further analysis.

    \begin{table}[ht!]
      \centering
      \footnotesize
      \begin{tabular}{r | r r | r r | r r}
          \toprule \multirow{2}{0.8cm}{\acp{dof}} & \multicolumn{2}{c}{No mask} & \multicolumn{2}{c}{Extracellular mask} & \multicolumn{2}{c}{Intracellular mask} \\ 
                   & Iters & Time    & Iters &   Time  & Iters & Time \\ \midrule
          9826     & 158   & {\bf 7.93}    & --    & --      & --    & --    \\
          71874    & 71    & {\bf 18.27}  & --    & --      & --    & --    \\
          549250   & 38    & {\bf 66.32}   & 50    & 83.07   & 45    & 74.55   \\
          4293378  & 46    & 626.35  & 40    & {\bf 521.78}  & 142   & 1686.44 \\ 
          15290746 & 45    & 2207.21 & 47    & {\bf 2184.99} & --   & -- \\ \bottomrule
      \end{tabular}
      \caption{Bidomain equations, physics-based mask test: Iteration count and CPU times in seconds for each of the masks considered (no mask, extracellular, intracellular) for increasing problem size. (--) denotes non-convergence. The lowest CPU time for each row has been highlighted with bold font.}\label{table:bidomain-mask}
    \end{table}

    \subsubsection{Bidomain: Adaptivity test}
    We show the iteration count and CPU times in seconds obtained when solving the Bidomain equations with all the proposed adaptivity strategies in Table~\ref{table:bidomain-adapt}. In this case, all adaptivity strategies yield a robust behavior that differs from the non-adapted case by very few iterations in all but the smallest problem. This confirms, as in the previous tests, that the adaptive masks are stable by construction.

    The difference in performance between both subselection strategies and the non-adapted case is not significant. This happens because the iteration counts are identical for all but the smallest problem, and there are many nonlinear terms in the residual, which makes assembly time more significant than the computation of the QR factorization. The randomized adaptivity approaches yields more significant improvements, as they require less iterations for convergence. In particular, the Randomized Constant strategy yields CPU time savings of up to a 17\% with respect to the non-adapted method in the largest problem.
    \begin{table}
        \centering
        \scriptsize
      \begin{tabular}{r | r r | r r | r r | r r | r r}
          \toprule \multirow{2}{0.8cm}{\acp{dof}}  & \multicolumn{2}{c|}{No adapt} & \multicolumn{2}{c|}{Sub pow} & \multicolumn{2}{c|}{Sub const} & \multicolumn{2}{c|}{Rand pow} & \multicolumn{2}{c}{Rand const} \\ 
                   & Iters & Time    & Iters & Time    & Iters  & Time    & Iters & Time    & Iters & Time    \\ \midrule
          9826     & 158   & {\bf 7.93}    & 269   & 12.34   & 269    & 12.37   & 306   & 14.01   & 393   & 18.31   \\
          71874    & 71    & 18.27   & 71    & {\bf 16.41}   & 71     & 16.69   & 74    & 16.93   & 76    & 17.71   \\
          549250   & 38    & 66.32   & 38    & 66.41   & 38     & 66.58   & 39    & {\bf 63.16}   & 39    & 63.36   \\
          4293378  & 46    & 626.35  & 46    & 633.10  & 46     & 631.26  & 47    & 605.57  & 47    & {\bf 596.32}  \\ 
          15290746 & 45    & 2207.21 & 45    & 2229.48 & 45     & 2219.03 & 44    & 2022.23 & 40    & {\bf 1838.64} \\ \bottomrule
      \end{tabular}

      \caption{Bidomain, adaptivity test: Iteration count and CPU times in seconds for increasing problem size and all adaptivity strategies considered (none, Subselect power, Subselect constant, Randomized power, and Randomized constant). The lowest CPU time for each row has been highlighted with bold font. }\label{table:bidomain-adapt}
      \end{table}

      \subsection{Bidomain: Best overall performance test}

      For this section, we solved the Bidomain equations with all combinations of physics-based masks, adaptivity strategies, and alternation ($p$ in $\{1,2,3\}$). The value $p=4$ resulted in divergence for most methods. For each problem size, we report the performance of \ac{aap}(50,1) and compare it with the best performing method in terms of CPU time. For this problem, it is consistently better to use alternation with $p=3$. In addition, the best strategy always involves adaptivity, and for the three largest problems this is always randomized adaptivity. For the largest problem this results in an iteration and CPU time reduction of 40\% and 43\% respectively. 

    \begin{table}[ht!]
      \centering
      \footnotesize
      \begin{tabular}{r | rr | r r | l l c }
          \toprule \multirow{2}{0.8cm}{\acp{dof}} & \multicolumn{2}{c|}{No mask} & \multicolumn{5}{c}{Best}  \\ 
                   & Iters & Time    &  Iters & Time    & Mask & Adaptivity & Alt ($p$)  \\ \midrule
          9826     & 158   & 6.13    & 123    & {\bf 5.27}    & No mask & Subselect power & $3$   \\
          71874    & 71    & 18.27   & 69     & {\bf 15.30}   & No mask & Subselect power & $3$ \\
          549250   & 38    & 66.32   & 33     & {\bf 53.75}   & No mask & Random power & $3$ \\
          4293378  & 46    & 631.26  & 27     & {\bf 343.26}  & Extracellular mask & Random constant & $3$ \\
          15290746 & 45    & 2207.21 & 27     & {\bf 1247.51} & No mask & Random power & $3$   \\ \bottomrule
      \end{tabular}
      \caption{Bidomain, best overall test: Iteration count and CPU time in seconds for the reference \ac{aa} method and for the best performing one. Best performing methods are characterized by their mask, adaptivity strategy, and alternation parameter $p$. The lowest CPU time for each row has been highlighted with bold font.}\label{table:bidomain-best}
    \end{table}

\subsection{Navier-Stokes}
    In this section we solve the stationary Navier-Stokes equation for the lid driven cavity benchmark: Given a Reynolds number $Re>0$, find a velocity $u:\Omega \to \R^2$ and a pressure $\pressureNS:\Omega \to \R$ such that
    \begin{equation}
        \begin{aligned}
            -\frac{1}{Re}\Delta \vec u + [\grad \vec u]\vec u + \grad \pressureNS &= 0  && \text{ in $\Omega$}, \\
            \dive \vec u &= 0 && \text{ in $\Omega$}, \\
            \vec u &= \vec u_D && \text{ on $\partial\Omega$}, \\
            \int_\Omega \pressureNS\,dx &= 0. 
        \end{aligned}
    \end{equation}
    We discretize this problem using inf-sup stable Taylor-Hood elements, and solve it using an accelerated Picard iteration as in \cite{pollock2019anderson}: Consider $Re>0$ and a given velocity $\vec u^k:\Omega\to\R$, then find a velocity $\vec u^{k+1}$ and a pressure $\pressureNS^{k+1}$ such that 
     \begin{equation}
        \begin{aligned}
            -\frac{1}{Re}\Delta \vec u^{k+1} + [\grad \vec u^{k+1}]\vec u^k + \grad \pressureNS^{k+1} &= 0  && \text{ in $\Omega$}, \\
            \dive \vec u^{k+1} &= 0 && \text{ in $\Omega$}, \\
            \vec u^{k+1} &= \vec u_D && \text{ on $\partial\Omega$}, \\
            \int_\Omega \pressureNS^{k+1}\,dx &= 0. 
        \end{aligned}
    \end{equation}
    The solution map $(\vec u^{k+1}, \pressureNS^{k+1}) = T(\vec u^k)$ yields a fixed point iteration whose fixed point is the solution of the stationary Navier-Stokes equation. This map $T$ is the one we accelerate with \ac{aap}. For this problem, we consider a physical mask study, an adaptivity study, and an overall best strategy computation. In all cases we considered $Re=5000$, a div-div stabilization with constant $\gamma=1$, and a direct solver for the linear system required at each Picard iteration. Convergence was established using a relative tolerance of $10^{-8}$.

    \subsubsection{Navier-Stokes: Physics-based mask test}
    In this test, we consider the same restriction operators $\mat\Pi_1$ from Stokes as defined in \eqref{eq:up-masks}, and show the iteration count and CPU time in seconds in Table~\ref{table:ns-mask} for all mask types. For this problem we observe that, in contrast to the Stokes problem, the pressure mask yields the worst performance, and the velocity mask the best ones. When comparing the velocity mask with the case without a mask, we observe that the difference in performance is negligible. This is the case because in this problem the dominant computational time is the direct solver required at each Picard iteration.

    It is important to highlight the difference between the fixed-point mapping used in this problem and the one from Stokes. In the Stokes problem, we considered a monolithic solver that can be recast as a pressure problem with the Schur complement. In this problem, our fixed point iteration is given by the map $\vec u^k \to (\vec u^{k+1}, \pressureNS^{k+1})$, where the problem can instead be recast in terms of the velocity by a simple restriction through $\vec u^k \to (\vec u^{k+1}, \pressureNS^{k+1}) \to \vec u^{k+1}$. This explains why for Stokes the pressure mask performed better than the velocity mask, whereas the velocity mask performs better than the pressure mask for Navier-Stokes.

    \begin{table}[ht!]
      \centering
      \footnotesize
      \begin{tabular}{r | r r | r r | r r}
          \toprule \multirow{2}{0.8cm}{\acp{dof}} & \multicolumn{2}{c|}{No mask} & \multicolumn{2}{c|}{Velocity mask} & \multicolumn{2}{c}{Pressure mask} \\ 
                  & Iters  &  Time    & Iters  &    Time & Iters  & Time  \\ \midrule
          9026    & 26 & 2.81   & 26  & {\bf 2.78}   & 33  & 3.48   \\
          36482   & 21 & 12.48  & 22  & 12.91  & 30  & {\bf 9.62}   \\
          146690  & 20 & 65.39  & 20  & {\bf 64.79}  & 24  & 76.74  \\
          588290  & 19 & 327.59 & 19  & {\bf 312.17} & 25  & 391.82 \\ 
          1324802 & 19 & {\bf 979.72} & 19  & 984.40 & 28 & 1423.89 \\ \bottomrule
      \end{tabular}
      \caption{Navier-Stokes, physics-based mask test: Iteration count and CPU times in seconds with no mask, velocity mask, and pressure mask, for increasing problem sizes. The lowest CPU time for each row has been highlighted with bold font.}\label{table:ns-mask}
    \end{table}

    \subsubsection{Navier-Stokes: Adaptivity test}
    We show the iteration count and CPU times in seconds required to solve the problem in Table~\ref{table:ns-adapt} for all adaptivity strategies considered. All methods behave similarly in terms of iteration count and verify the stability of the adaptive strategies. The mild increase in iterations related to the subselection strategies results in worse CPU times. Instead, the randomized strategies result in the same number of iterations and similar time-to-solution as the non-adapted case. Differences in CPU times for this test are negligible for methods with similar iteration counts, and the cases in which these differences are more significant are mainly due to the number of iterations. This happens because the dominant cost for this method is the computation of the LU factorization required to solve the linearized problem at each Picard iteration. 

    \begin{table}
        \centering
        \footnotesize
      \begin{tabular}{r | r r | r r | r r | r r | r r}
          \toprule \multirow{2}{0.8cm}{\acp{dof}} & \multicolumn{2}{c|}{No adapt} & \multicolumn{2}{c|}{Sub pow} & \multicolumn{2}{c|}{Sub const} & \multicolumn{2}{c|}{Rand pow} & \multicolumn{2}{c}{Rand const} \\ 
                 & Iters &  Time  & Iters & Time    & Iters & Time    & Iters & Time   & Iters & Time  \\ \midrule
          9026    & 26    & 2.81   & 28    & 2.99    & 28    & 2.98    & 27    & 2.88   & 26    & {\bf 2.78}  \\
          36482   & 21    & {\bf 12.48}  & 22    & 12.91   & 22    & 12.92   & 21    & 12.37  & 22    & 12.92 \\
          146690  & 20    & 65.39  & 21    & 68.03   & 21    & 68.85   & 20    & {\bf 59.95}  & 20    & 65.03 \\
          588290  & 19    & {\bf 327.59} & 25    & 423.38  & 25    & 422.96  & 19    & 328.60 & 19    & 328.48 \\
          1324802 & 19    & {\bf 979.72} & 25    & 1276.82 & 25    & 1273.59 & 19    & 997.64 & 19    & 988.47 \\  \bottomrule
      \end{tabular}

      \caption{Navier-Stokes, adaptivity test: Iteration count and CPU times in seconds for all considered adaptivity strategies (None, Subselect power, Subselect constant, Randomized power, and Randomized constant) and for increasing problem size. (--) denotes non-convergence. The lowest CPU time for each row has been highlighted with bold font.}\label{table:ns-adapt}
      \end{table}

      \subsection{Navier-Stokes: Best overall performance test}

      In this section, we tested all combinations of physics-based masks (none, velocity, pressure), adaptivity strategies, and alternation parameter ($p$ in $\{1,2,3,4\}$). For each problem size, we report in Table~\ref{table:ns-best} the performance of the base method and compare it with the best performing one in terms of CPU time. In this test, we see that it is always convenient to consider a velocity mask together with either adaptivity, alternation, or both. This results in a reduction of 1-2 iterations, which reduces the computational cost by up to a mild 9\%.

    \begin{table}[ht!]
      \centering
      \footnotesize
      \begin{tabular}{r | rr | r r | ll c }
          \toprule \multirow{2}{0.8cm}{\acp{dof}} & \multicolumn{2}{c|}{No mask} & \multicolumn{5}{c}{Best}  \\ 
                  & Iters &  Time  & Iters & Time   & Mask & Adaptivity & Alt ($p$)  \\ \midrule
          9026    &  26   & 2.81   & 26    & {\bf 2.78}   & Velocity mask & No adaptivity & $1$   \\
          36482   &  21   & 12.48  & 20    & {\bf 11.67}  & Velocity mask & No adaptivity & $2$  \\
          146690  &  20   & 65.39  & 18    & {\bf 54.10}  & Velocity mask & Subselect power & $1$   \\
          588290  &  19   & 327.59 & 18    & {\bf 291.60} & Velocity mask & No adaptivity & $2$         \\
          1324802 &  19   & 979.72 & 18    & {\bf 892.88} & Velocity mask & Subselect power & $2$ \\  \bottomrule
      \end{tabular}
      \caption{Navier-Stokes, best overall test: Iteration count and CPU time in seconds for \ac{aa} and the best performing method among the ones considered. The best methods are characterized by their mask, adaptivity strategy, and alternation parameter $p$. The lowest CPU time for each row has been highlighted with bold font.}\label{table:ns-best}
    \end{table}

\section{Conclusions}\label{section:conclusions}
In this work we proposed a two-level \ac{aap} algorithm. The two-level nature of the method is motivated by the possibility of using two separate restriction operators, where one of them is physics-based and known before execution, and the other one is purely algebraic, adaptive, and dynamically computed at each iteration. This first operator is advantageous as it reduces the memory requirements of the method, and in some cases it was observed to yield a significant reduction in the number of iterations.  The second operator further reduces the size of the least-squares problem by adaptively choosing \acp{dof}. The adaptivity strategy that we proposed was directly motivated by the backward stability analysis, and we verified in all tests that the strategies were in fact stable, as the difference in iterations between the method with and without adaptivity was typically not very significant, and proved beneficial in most cases. In addition, the Stokes and Navier-Stokes tests reveal that the physics-based mask considered depends heavily on the structure of the fixed-point operator that one accelerates, as both models are based on similar variables but were solved with very different strategies.

{We emphasize that the purpose of the present analysis is not to establish global convergence rates for Anderson acceleration under sketching, but to characterize admissible approximations that preserve the convergence behavior of the underlying exact method. A general convergence theory for nonlinear Anderson acceleration remains an open problem, even in the absence of sketching.}

The proposed algorithm was shown to be effective for certain problems, where we observed time savings between a $40\%$ and a $50\%$ with respect to standard \ac{aa}. These savings are not consistent among all tests, which shows that performance is problem dependent, as one would expect. The problems where our algorithm was most advantageous are the Stokes and the Bidomain equations, which are the ones where residual evaluation was the cheapest. This shows that our methodology can improve the performance of \ac{aa} most reliably in problems where the computation of the residual is not excessively dominant, as otherwise the least-squares problems solved in \ac{aap} provides no significant overhead, and a minimal increase in iterations can be disadvantageous. Despite this, in all cases we could obtain a method that was comparable or better than the classic one in terms of iterations, and it remained robust with respect to the problem size despite the reduction of the least-square problem.

{Although the observed speedups are moderate, such reductions still significantly accelerate the convergence of the numerical scheme in large-scale PDE simulations where AA is applied repeatedly and the least-squares solve constitutes a non-negligible fraction of the total runtime. Moreover, the proposed approach achieves these gains without altering the underlying nonlinear solver or preconditioner, making it complementary to existing algorithmic improvements.}

Future work will integrate our Julia package with more advanced adaptive sketching techniques provided in RandBLAS, the Randomized Numerical Linear Algebra (RNLA) package in BLAS \cite{murray2023randomized} currently under development.  This work has been implemented in Julia, and will be released in the AAP.jl library, available to the community upon the publication of this work. 

\section{Acknowledgments}
NAB was supported by Centro de Modelamiento Matematico (CMM), Proyecto Basal FB210005, by ANID Postdoctoral Proyecto 3230326, and by the SIAM Postdoctoral grant. MLP was sponsored by the Artificial Intelligence Initiative as part of the Laboratory Directed Research and Development (LDRD) Program of Oak Ridge National Laboratory, managed by UT-Battelle, LLC, for the US Department of Energy under contract DE-AC05-00OR22725. 

\bibliography{main}
\bibliographystyle{alpha}

\appendix
\section{Sensitivity on the summable sequence}\label{appendix:power}

In Theorem \ref{thm:adaptive}, we obtained that backward stability depends on an arbitrary summable sequence $(\eta_j)_j$, for which we have chosen $\eta_j = j^{-\ell}$ with $\ell=1.1$. One might naturally wonder whether this arbitrary exponent is fundamental for the obtained results. To answer this question, we chose the Bidomain problem with 71874 \acp{dof} and solved it with varying values of $\ell$ using both subselect and randomized masks. We show the results in Table \ref{tab:power-decay} where we report the total number of iterations and the number of adaptive steps for both subselection and randomized masks. We observe that (1) higher exponents effectively deactivate adaptivity, (2) the subselection mask is robust with respect to this parameter, and that (3) randomized masks behave better when we stay within the theory $\ell > 1$, with only a minor deterioration otherwise. We conclude that it is safe to use $\ell=1.1$ as it balances theoretical correctness while also enabling adaptivity in practice.
    \begin{table}[ht!]
      \centering
      \footnotesize
      \begin{tabular}{r | cc | cc}
          \toprule $\ell$ & \multicolumn{2}{c|}{Subselection mask} & \multicolumn{2}{c}{Randomized mask}  \\ 
               & Iterations &  Adaptive steps  & Iterations & Adaptive steps \\ \midrule
          0    & 71 & 50 &  81 & 78  \\
          1    & 71 & 50 &  76 & 73  \\
          1.1  & 71 & 50 &  76 & 73  \\
          2    & 71 & 49 &  77 & 74  \\
          5    & 71 & 0  &  71 & 16  \\
          10   & 71 & 0  &  71 & 1  \\
          100  & 71 & 0  &  71 & 0  \\
          1000 & 71 & 0  &  71 & 0  \\  \bottomrule
      \end{tabular}
      \caption{Sensitivity study on the coefficient $\ell$. The first two rows ($\ell$ in $\{0,1\}$) are outside of the theory, and all others yield a summable sequence. For each mask type (subselection and randomized), we provide both the total iterations and the number of iterations in which adaptivity was used.}\label{tab:power-decay}
    \end{table}

\end{document}